\documentclass{amsart}
\usepackage{amsmath}
\usepackage{amscd}
\usepackage{amssymb}
\usepackage{graphicx}
\usepackage[colorlinks=true]{hyperref}

\newtheorem{theorem}{Theorem}

\newtheorem{lemma}[theorem]{Lemma}
\newtheorem{proposition}[theorem]{Proposition}
\newtheorem{definition}[theorem]{Definition}

\newcommand{\diam}{\operatorname{diam}}

\newcommand{\R}{\mathbb{R}}

\newcommand{\N}{\mathbb{N}}
\newcommand{\G}{\Gamma}

\title[Remarks on quasi-isometric non-embeddability]
{\large Remarks on quasi-isometric non-embeddability into
uniformly convex Banach spaces$^{\dagger}$}

\author{Piotr W. Nowak}

\thanks{$^\dagger$ This research was partially supported by a Summer Research Grant from the College of Arts \& Science of Vanderbilt University}
\address{Department of Mathematics, Vanderbilt University, 1326 Stevenson Center, Nashville, TN 37240
USA.}

\email{pnowak@math.vanderbilt.edu}


\begin{document}

\begin{abstract}
We construct a locally finite graph and a bounded geometry metric
space which do not admit a quasi-isometric embedding into any
uniformly convex Banach space. Connections with the geometry of
$c_0$ and superreflexivity are discussed.
\end{abstract}
 \maketitle
The question of coarse embeddability into uniformly convex Banach
spaces became interesting after the recent work of G.~Kasparov and
G.~Yu, who showed the coarse Novikov Conjecture, i.e. that the
coarse assembly map in K-theory is injective, for bounded geometry
metric spaces coarsely embeddable into uniformly convex Banach
spaces \cite{kasparov-yu}. So far there is no example of a bounded
geometry metric space which wouldn't admit a coarse embedding into
any uniformly convex Banach space - such spaces are hard to find
even if we restrict the target space to be Hilbert.

In this note we consider quasi-isometric embeddings, a special
case of coarse embedding that may be described as \emph{large
scale biLipschitz}. We construct a locally finite graph which does
not admit a quasi-isometric embedding into any uniformly convex
Banach space. In contrast to this note that every separable metric
space admits a biLipschitz embedding into a strictly convex Banach
space, namely into the space $c_0$ with an equivalent strictly
convex norm, by a classic result of I.~Aharoni \cite{aharoni}.

It also turns out that our methods can be applied to a $c_0$-type
of geometry. As a result we find an explicit geometric obstruction
to uniform quasi-isometric embeddability (i.e. with embedding
constants independent of $n$) of the $\ell_{\infty}^n$'s into any
uniformly convex Banach space. In the proof one can directly see
how the Lipschitz constants of the embedding are being pushed away
to infinity as $n$ grows larger.

In the last section we also comment our observations in view of
Bourgain's paper on superreflexivity \cite{bourgain}, from which
quasi-isometric non-embeddability of trees into uniformly convex
Banach spaces can be deduced. Bourgain's proof is of different
nature, using classic characterization of superreflexivity due to
R.C.~James and G.~Pisier. We show that our graph has essentially
the opposite of the hyperbolic geometry on the tree and that the
two results are independent of each other.

The inspiration for the construction of the graph $\Gamma$ comes
from the work of T.~Laakso on $A_{\infty}$-weighted metrics on the
plane \cite{laakso}, similar spaces were also considered in
\cite{newman-rabinovich} in the context of metric embeddings of
finite graphs into Hilbert spaces. I would like to greatly thank
Piotr Haj\l asz from whom I learned about Laakso's work and my
advisor Guoliang Yu for  his support and guidance. I am also very
grateful to William Johnson for illuminating discussions in which
the connection with Banach space geometry (among other things)
became apparent.

\section{Construction of the graph $\G$}
The main idea is to build a graph with a fractal-like structure in
such way that the deformed geometry is reflected on the large scale.
We will first inductively construct a sequence
$\left\lbrace\G\right\rbrace_{n\in \N}$ of finite graphs which will
be the building blocks of $\G$ (compare \cite{newman-rabinovich}).\\

Define $\G_0$ to be a single edge of length 1. To construct $\G_1$
take four copies of $\G_0$ and denote by $p_i$ and $q_i$ the
vertices of $i$-th copy, $i=0,1,2,3$. The graph $\G_1$ is
constructed by identifying $q_i$ with $p_{i+1}$ with $i \mod 4$.
Equip $\G_1$ with the path metric.

Similarly to construct $\G_{n+1}$ take four copies of $\G_n$,
denote two vertices of valence $2^{n}$ that are distance $2^{n}$
away from each other in the $i$-th copy by $p_i$ and $q_i$,
$i=0,1,2,3$ (there are two pairs of such vertices in $\G_n$), and
identify $q_i$ with $p_{i+1}$ with $i \mod 4$. Equip each $\G_n$
with a path metric. For each $\G_n$ a pair of vertices of valence
$2^n$ and distance $2^n$ from each other will be called \emph{a
pair of primary vertices}. In every $\G_n$ there are exactly two
pairs of primary vertices.
Note that $\diam\G_n=2^n$.\\

For each $n\in\N$ denote by $p_n$ and $q_n$ a pair of primary
vertices in $\G_n$. Construct the graph $\G$ by identifying $q_n$
with $p_{n+1}$ for every $n\in \N$  (see Fig. 2) and extending the
metric in the obvious way. Note that one of the pairs of primary
vertices in $\G_n$ changes valence under the isometric embeddings
$\G_n\subseteq \G_{n+1}$ and $\G_n\subseteq \G$, however we will
refer to them without change as primary vertices of an isometric
copy of $\G_n$.

\begin{figure}
\begin{center}
\includegraphics{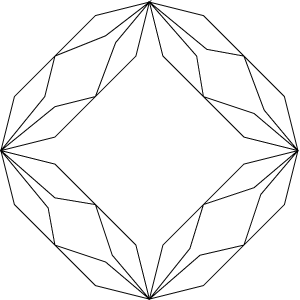}\label{fig:gammathree} \caption{ The
graph $\G_3$}
\end{center}
\end{figure}

Actually, for our purposes it would be enough to consider just the
set of all vertices of $\Gamma$, however we found the quote from
the Introduction to \cite{gromov}: \emph{`Given a discrete metric
space $\G$, one can make it more palatable by adding some meat to
$\G$ in the form of edges and higher dimensional simplices with
vertices in $\G$, without changing the quasi-isometry type'},
quite applicable in this situation.\\

Recall that a metric space $X$ is called locally finite if there
exists a discrete subset $\mathcal{N}\subseteq X$ and a constant
$C>0$ such that for every $x\in X$ there is a $y\in \mathcal{N}$
satisfying $d(x,y)\le C$ and for every $y\in \mathcal{N}$ the number
of elements in every ball around $y$ in $\mathcal{N}$ is finite. The
graph $\G$ is a
locally finite metric space.\\

\begin{figure}
\begin{center}
\includegraphics[width=5in]{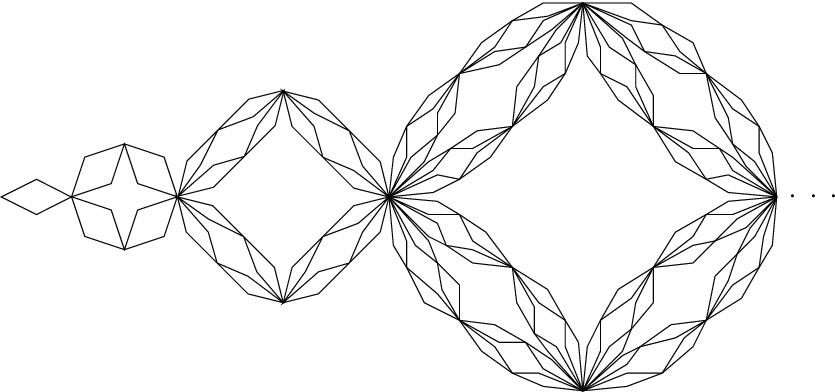}\label{fig:gammawhole}\caption{ The
graph $\G$}
\end{center}
\end{figure}

\section{Round ball spaces and nonexistence of the embedding}
In this section we prove that the graph $\G$ constructed above does
not admit a quasi-isometric embedding into any \emph{round ball
metric space}, which we define below.

\begin{definition}\label{def.qi}
\normalfont Let $X,Y$ be metric spaces. We say that $f\colon X\to Y$
is a \emph{quasi-isometry} if there are constants $L>0$, $C\ge 0$
such that
$$L^{-1}d_X(x,y) -C\le d_Y(f(x),f(y))\le L d_X(x,y) +C$$
for all $x,y\in X$.\\
\end{definition}

\begin{definition}[\cite{laakso}]\normalfont
A metric space $X$ is called a \emph{round ball space} if for every
$\varepsilon>0$ there exists $\delta_{\varepsilon}>0$ such that
$$\diam\left(B\left(x, \frac{1+\delta_{\varepsilon}}{2}d(x,y)\right)\cap B\left(y, \frac{1+\delta_{\varepsilon}}{2}d(x,y)\right)
\right)<\varepsilon d(x,y)$$ for all $x,y\in X$.\\
\end{definition}

The round ball condition generalizes to metric spaces the classical
notion of uniform convexity for Banach spaces. This is made precise
in the following

\begin{proposition}[\cite{laakso}]
A Banach space is a round ball space if and only if it is
uniformly convex.\\
\end{proposition}

Given a map $f:X\to Y$ between metric spaces and points $x,y\in X$,
by $L_{x,y}$ we denote the Lipschitz constant of $f\vert_{\{x,y\}}$,
i.e.
$$L_{x,y}=\frac{d_Y(f(x),f(y))}{d_X(x,y)}.$$
The following lemma can be extracted from \cite{laakso}.

\begin{lemma}\label{lemma.of.4.points}
Let $x_1,x_2, x_3, x_4\in X$ satisfy
\begin{enumerate}
\item $d_X(x_1,x_3)=d_X(x_2,x_4)=C$\label{condition.1.lemma}
\item
$d_X(x_1,x_2)=d_X(x_1,x_3)=d_X(x_2,x_4)=d_X(x_3,x_4)=C/2$\label{condition.2.lemma}
\end{enumerate}
for some constant $C>0$. If $f:X\to Y$ is a map to a round ball
metric space $Y$ and the restriction $f\vert_{\left\lbrace
x_1,x_2,x_3,x_4\right\rbrace}$ is biLipschitz with constant $L$ then
$$\max\left\lbrace L_{x_1, x_2}, L_{x_1, x_4}, L_{x_3, x_2}, L_{x_3, x_4}\right\rbrace \ge (1+\delta_{L^{-2}})
L_{x_1,x_3}.$$
\end{lemma}
\begin{proof}
Condition (\ref{condition.1.lemma}) implies that $$d_Y(f(x_2),
f(x_4))\ge L^{-2} d_Y(f(x_1), f(x_3)).$$ If $B(f(x_1), R)\cap
B(f(x_3, R)$ contains $f(x_2)$ and $f(x_4)$ then we must have $R\ge
\frac{1}{2}\left(1+\delta_{L^{-2}}\right)d_Y(f(x_1), f(x_3))$, by
the round ball condition. Since this is the case when we take $$R =
\max\left\lbrace d_Y(f(x_1),f(x_2)), d_Y(f(x_1),f(x_4)),
d_Y(f(x_3),f(x_2)), d_Y(f(x_3),f(x_4))\right\rbrace,$$ the assertion
follows.
\end{proof}

\begin{theorem}\label{theorem.nonembed}
The space $\G$ does not admit a quasi-isometric embedding into any
round ball space.
\end{theorem}

\begin{proof}
Assume that there exists a quasi-isometric embedding of the metric
space $\G$ into a round ball space $Y$. Observe that such an
embedding is \emph{biLipschitz for large distances}, i.e. there
are constants $L>0$ and $S>0$ such that $$L^{-1} d_{\G}(x,y)\le
d_Y(f(x),f(y))\le L d_{\G}(x,y)$$ for all $x,y\in \G$ satisfying
$d_{\G}(x,y)\ge S$. Choose $n\in \N$ such that
$\left(1+\delta_{L^{-2}}\right)^nL^{-1}> L$. In $\G$ choose points
$x_1^{(0)}$ and $x_3^{(0)}$ that are images of primary vertices in
an isometric copy of $\G_k$ for some $k\in\N$ and satisfy
$d_{\G}(x_1^{(0)},x_3^{(0)})> 2^{n+1} S$. Denote by $x_2^{(0)}$
and $x_4^{(0)}$ the second pair of primary vertices in this copy
of $\G_k$. The quadruple of points $x_1^{(0)}$, $x_2^{(0)}$,
$x_3^{(0)}$, $x_4^{(0)}$ satisfies the hypothesis of lemma
\ref{lemma.of.4.points} and we get
\begin{equation}\label{equation.max.lip.const}
\max\left\lbrace L_{x_1^{(0)},x_2^{(0)}}, L_{x_1^{(0)},x_4^{(0)}},
L_{x_3^{(0)},x_2^{(0)}}, L_{x_3^{(0)},
x_4^{(0)}}\right\rbrace>(1+\delta_{L^{-2}}) L_{x_1^{(0)}, x_3^{(0)}}
\end{equation}
Rename the pair for which the $\max$ on the left side is attained to
$x_1^{(1)}$ and $x_3^{(1)}$. These two points constitute a pair of
primary vertices of an isometric copy of $\G_{k-1}$. Denote by
$x_2^{(1)}$ and $x_4^{(1)}$ the remaining pair of primary vertices
in $\G_{k-1}$ . Applying lemma \ref{lemma.of.4.points} to the
quadruple $x_1^{(1)}$, $x_2^{(1)}$, $x_3^{(1)}$, $x_4^{(1)}$ and the
renaming procedure together with equation
(\ref{equation.max.lip.const}) we get a pair of points $x_1^{(2)}$
and $x_3^{(2)}$ satisfying
$$L_{x_1^{(2)},x_3^{(2)}}\ge (1+\delta_{L^{-2}})^2 L_{x_1^{(0)},
x_3^{(0)}}.$$ The points $x_1^{(2)}$ and $x_3^{(2)}$  are again a
pair of primary vertices of an isometric copy of $\G_{k-2}$.

Continuing in this way after $n$ steps we will get a pair of points
$x_1^{(n)}$, $x_3^{(n)}$ satisfying $d_{\G}(x_1^{(n)}, x_3^{(n)})\ge
S$ and
$$L_{x_1^{(n)},x_3^{(n)}}\ge (1+\delta_{L^{-2}})^n L_{x_1^{(0)},
x_3^{(0)}}.$$ By the choice of $n$ we get a contradiction to the
Lipschitz condition for large distances.
\end{proof}

We remark that what is essential for the proof of Theorem
\ref{theorem.nonembed} is the existence of an isometric copy of
$\G_n$ in $\G$ for arbitrarily large $n$. Given the sequence
$\left\lbrace \G_n\right\rbrace_{n\in\N}$ one can thus create
similar examples using constructions like e.g. disjoint union with
an
appropriate metric.\\\\
\emph{Bounded geometry example.} We also want to indicate that it
is now easy to modify the construction of $\G$ to get a bounded
geometry metric space with the property of non-embeddability.
Simply consider for each $n\in\N$ the space $V_n$, the set of
vertices of $\G_n$ with the metric multiplied by $n$ and glue them
together as in the construction of $\G$. The resulting space
however is not
quasi-geodesic.\\

\section{The geometry of $c_0$ and quasi-isometric embeddings}
The fact that a 1-net in $c_0$ does not admit a biLipschitz for
large distances embedding into any uniformly convex Banach space
is well known in Banach space geometry, one can prove it using
e.g. ultrapowers. We will show however how the methods from the
previous section can be implemented in the $\ell_{\infty}^n$'s and
$c_0$ and we will give an explicit geometric obstruction to
uniform quasi-isometric embeddability (i.e. with constants $L$ and
$C$ independent of $n$)\footnote{We do not use here the standard
notion of uniform containment of a sequence $X_n$ of finite
dimensional spaces in a Banach space, namely that for every
$\epsilon>0$ there is a isomorphic embedding of $X_n$ with
distortion less than $1+\epsilon$ for every $n$, since it
emphasizes the infinitesimal aspect of uniformity. Our definition
is suitable for the purposes of large scale geometric behavior.}
of $\ell_{\infty}^n$'s into uniformly convex Banach spaces.

Consider the embedding of $\ell_{\infty}^n$ into
$\ell_{\infty}^{n+1}$ given by adding $(n+1)$-st coordinate 0 to
each vector in $\ell_{\infty}^n$. For vectors
$v,w\in\ell_{\infty}^n$ take their images under this embedding
$\widetilde{v}, \widetilde{w}\in \ell_{\infty}^{n+1}$ and the
vectors $x=\frac{\widetilde{v}+\widetilde{w}}{2}+
\big(\underbrace{0,0,...,0}_n,\frac{\Vert
\widetilde{v}-\widetilde{w}\Vert}{2}\big)$,
$y=\frac{\widetilde{v}+\widetilde{w}}{2}+
\big(\underbrace{0,0,...,0}_n,-\frac{\Vert
\widetilde{v}-\widetilde{w}\Vert}{2}\big)$. It is easy to check
that the quadruple of points $\widetilde{v},\widetilde{w},x,y$
satisfies the hypothesis of lemma \ref{lemma.of.4.points}, thus we
can apply with no change the procedure from the proof of Theorem
\ref{theorem.nonembed} and recover

\begin{proposition}
Let $X$ be a Banach space containing $\ell_{\infty}^n$'s
quasi-isometrically uniformly. Then $X$ does not admit a
quasi-isometric embedding into any uniformly convex Banach space.
\end{proposition}

The same argument gives a purely metric proof of the fact that the
unit ball in $c_0$ does not admit a biLipschitz embedding into any
uniformly convex Banach space, which is again obvious once we
appeal to the linear structure of $c_0$.

\section{Some remarks on a paper of Bourgain}

For the purposes of the coarse Novikov Conjecture
\cite{kasparov-yu} one can consider coarse embeddings into
superreflexive Banach spaces, since these are exactly the ones
that admit an equivalent uniformly convex norm, due to a theorem
of Enflo \cite{enflo}. It might be thus interesting to confront
the above observations with a paper of J.~Bourgain in which a
metric characterization of superreflexivity was given
\cite{bourgain}. The necessary and sufficient condition for a
Banach space to be superreflexive was shown to be biLipschitz
uniform non-embeddability of a sequence of trees $T_j$ defined
below (see \cite{bourgain} for a precise formulation).

Denote $\Omega_n=\lbrace -1,1\rbrace^n$, $T_n=\bigcup_{i\le n}
\Omega_i$, $T=\bigcup_{j=1}^{\infty} T_j$ and again add "some
meat" to $T$ in the form of edges in the obvious way and denote
the resulting space by $\overline{T}$. First note that
quasi-isometric embeddability of $T$ (or equivalently
$\overline{T}$) into a superreflexive Banach space implies the
existence of a biLipschitz embedding of $T$ into such a space.

\begin{lemma}
Let $X$ be a discrete metric space and assume that $X$ embeds
quasi-isometrically into a superreflexive Banach space. Then $X$
admits a biLipschitz embedding into a superreflexive Banach space.
\end{lemma}

\begin{proof}
Let $f:X\to E$ be the quasi-isometric embedding.
Define\footnote{Although it is not of great importance, it is
convenient to take the direct sum with the $\ell_1$-norm}
$\widetilde{f}:X\to E\oplus \ell_2(X)$ by the formula
$\widetilde{f}(x)=f(x)\oplus\delta_x$. $E\oplus \ell_2(X)$ is
super-reflexive and it is easy to verify that $\widetilde{f}$ is a
biLipschitz embedding.
\end{proof}
Thus $\overline{T}$ also does not admit a quasi-isometric
embedding into any superreflexive Banach space, the argument in
\cite{bourgain} is however of probabilistic nature, as mentioned
earlier.

The geometry of the graph $\G$ is intuitively the very opposite of
the hyperbolic geometry on a tree. This is indicated already by
topological invariants, but also the geometries of these spaces
are very different on the large scale. The next proposition shows
that Theorem \ref{theorem.nonembed} and Bourgain's result cannot
be deduced from each other.

\begin{proposition}
\begin{enumerate}
\item $\overline{T}$ does not admit a coarse embedding into $\G$
\item $\G$ does not admit a coarse embedding into $\overline{T}$.
\end{enumerate}
\end{proposition}

\begin{proof}
To see (1) note the obvious fact that since the graph is infinite
in just one direction, a coarse embedding from $\R$ to $\G$ must
map both infinite ends of $\R$ in the same direction in $\G$ and
distant points on the line have to cross close to the joints of
$\G$, i.e. the points in which $\G_n$ is glued with $\G_{n+1}$.
Since $\R$ is isometrically embedded in the tree $\overline{T}$,
the assertion follows.

Similarly for (2) take two infinite geodesic rays in $\G$ and
observe that they have to pass through the joints so that
arbitrarily far on these geodesic rays some points are identified.
In an isometric copy of $\G_n\subseteq \G$ the pair of primary
vertices that are not joints is at distance $2^n$ apart. Thus a
coarse embedding of $\G$ into $\overline{T}$ from some point on
cannot map these vertices into the same branch in $\overline{T}$.
But that means that the joints as points on two geodesic rays will
be mapped to points whose distance grows to infinity, which gives
a contradiction.
\end{proof}
Note that to get quasi-isometric non-embeddability of $\G$ into
$\overline{T}$ it is enough to use Theorem \ref{theorem.nonembed},
since any $\R$-tree is a round ball space.\\

 The intuition behind the facts presented above might be the following: the tree $T$
carries an $\ell_1$-geometry, while the geometry of $\G$ resembles
the one of $c_0$. Recall that $\ell_1$ and $c_0$ are the standard
examples of non-uniformly convex Banach spaces, thus the presence
of any of these geometries should be an obstruction to
quasi-isometric embeddability into uniformly convex Banach spaces.

\bibliographystyle{amsalpha}

\begin{thebibliography}{11111111}
\baselineskip=12pt

\bibitem[Ah]{aharoni}
\textsc{I.~Aharoni}, \emph{Every separable metric space is
Lipschitz equivalent to a subset of $c_0$}, Isr. J. Math. 19
(1974), 284-291.

\bibitem[Bo]{bourgain}
\textsc{J.~Bourgain}, \emph{The metrical interpretation of
superreflexivity in Banach spaces}, Isr. J. Math. 56 (1986),
222-230.

\bibitem[En]{enflo}
\textsc{P.~Enflo}, \emph{Banach spaces which can be given an
equivalent uniformly convex norm}, Isr. J. Math. 13 (1972),
281-288.

\bibitem[Gr]{gromov}
\textsc{M.~Gromov}, \emph{Asymptotic invariants of infinite
groups}, London Mathematical Society Lecture Notes, no.182, s.
1-295, Cambridge University Press, 1993.

\bibitem[KY]{kasparov-yu}
\textsc{G.~Kasparov, G.~Yu}, \emph{The coarse geometric Novikov
conjecture and uniform convexity}, preprint 2004.

\bibitem[La]{laakso}
\textsc{T.~Laakso}, \emph{Plane with $A_{\infty}$-weighted metric
not bilipschitz embeddable into $\R^n$}, Bull. London Math. Soc. 34
(2002) 667-676.

\bibitem[NR]{newman-rabinovich}
\textsc{I.~Newman, Y.~Rabinovich}, \emph{A Lower Bound on the
Distortion of Embedding Planar Metrics into Euclidean space},
Discrete Comput. Geom. 29:77-81 (2003).
\end{thebibliography}

\end{document}